\documentclass{amsart}\usepackage{amssymb,euscript,amsfonts,bbm,graphicx}
\usepackage{enumerate}
\newtheorem{Lemma}{Lemma}[section]\newcommand{\bel}{\begin{Lemma}}\newcommand{\eel}{\end{Lemma}}
\newtheorem{Example}[Lemma]{Example}\newcommand{\bex}{\begin{Example}\rm}\newcommand{\eex}{\end{Example}}
\newtheorem{Proposition}[Lemma]{Proposition}\newcommand{\bprop}{\begin{Proposition}}\newcommand{\eprop}{\end{Proposition}}
\newtheorem{Definition-Proposition}[Lemma]{Definition-Proposition}

\def\bpr{~\\{\em Proof.\ }}\newcommand{\epr}{$\bull$\\}
\newtheorem{Theorem}[Lemma]{Theorem}\newcommand{\bthe}{\begin{Theorem}}\newcommand{\ethe}{\end{Theorem}}
\newtheorem{Definition}[Lemma]{Definition}\newcommand{\bed}{\begin{Definition}}\newcommand{\eed}{\end{Definition}}
\newtheorem{Remark}[Lemma]{Remark}\newcommand{\beR}{\begin{Remark}\rm}\newcommand{\eeR}{\end{Remark}}
\newtheorem{Remarks}[Lemma]{Remarks}\newcommand{\beRs}{\begin{Remarks}\rm}\newcommand{\eeRs}{\end{Remarks}}
\newtheorem{Corollary}[Lemma]{Corollary}\newcommand{\bcor}{\begin{Corollary}}\newcommand{\ecor}{\end{Corollary}}

\newcommand{\beq}{\begin{equation}}\newcommand{\eeq}{\end{equation}}
\newcommand{\beqn}{\begin{equation*}}\newcommand{\eeqn}{\end{equation*}}
\newcommand{\bem}{\begin{displaymath}}\newcommand{\eem}{\end{displaymath}}
\newcommand{\beqa}{\begin{eqnarray}}\newcommand{\eeqa}{\end{eqnarray}}
\newcommand{\bee}{\begin{enumerate}}\newcommand{\eee}{\end{enumerate}}
\newcommand{\bei}{\begin{itemize}}\newcommand{\eei}{\end{itemize}}
\newcommand{\bet}{\begin{tabular}{cccccccc}}\newcommand{\eet}{\end{tabular}}
\newcommand{\bpm}{\begin{pmatrix}}\newcommand{\epm}{\end{pmatrix}}
\newcommand{\bM}{\begin{matrix}}\newcommand{\eM}{\end{matrix}}
\newcommand{\ber}{\begin{array}{l}}\newcommand{\eer}{\end{array}}

\def\p{o}\def\o{o}\newcommand{\frD}{\mathfrak{D}}

\def\bull{\vrule height .9ex width .9ex depth -.1ex }
\newcommand{\quotient}[2]{{\left.\raisebox{1.6ex}{$#1$}\!\!\!\!\!{\scalebox{2}{\ensuremath\diagup}}
\!\!\!\!\!\raisebox{-1ex}{$#2$}\right.}}

\newcommand{\quotients}[2]{{\footnotesize\left.\raisebox{0.4ex}{$#1$}\! / \!\raisebox{-0.4ex}{$#2$}\right.}}

\def\di{\partial}
\def\ra{\rightarrow}

\def\into{\stackrel{i}{\hookrightarrow}}\def\proj{\stackrel{\pi}{\ra}}

\def\cN{\mathcal{N}}\def\cO{\mathcal{O}}
\def\cS{\mathcal{S}}
\def\cU{\mathcal{U}}

\def\cm{{\frak m}}

\def\A{\mathbb{A}}\def\C{\mathbb{C}}

\def\k{\mathbbm{k}}\def\P{\mathbb{P}}

\def\Z{\mathbb{Z}}

\def\De{\Delta}
\def\ep{\epsilon}\def\la{\lambda}
\def\Si{\Sigma}\def\Om{\Omega}

\def\tf{{\tilde{f}}}

  \def\tX{{\tilde{X}}}

\def\ux{\underline{x}}\def\uy{\underline{y}}

\def\empty{\varnothing}
\def\suml{\sum\limits}\def\oplusl{\mathop\oplus\limits}

\def\capl{\mathop\cap\limits}\def\cupl{\mathop\cup\limits}

\def\prodl{\prod\limits}\def\otimesl{\mathop\otimes\limits}

\def\lSi{{\overline{\Sigma}}}

\def\smin{\setminus}\def\sset{\subset}\def\sseteq{\subseteq}

\def\Pic{\rm{Pic}}
\def\omp{ordinary multiple point}\def\Db{{\De^\bot}}

\def\?{{\bf ???}}

\title{D\MakeLowercase{iscriminant of the ordinary transversal singularity type.}
T\MakeLowercase{he global equivalence class.}}
\author{M\MakeLowercase{axim} K\MakeLowercase{azarian}, D\MakeLowercase{mitry} K\MakeLowercase{erner} \MakeLowercase{and} A\MakeLowercase{ndr\'as} N\MakeLowercase{\'emethi}}
\date{\today}

\address{Steklov Mathematical Institute, RAS and the Poncelet Laboratory, Independent University of Moscow}
\email{kazarian@mccme.ru}
\address{Department of Mathematics, Ben Gurion University of the Negev, Israel}
\email{dmitry.kerner@gmail.com}
\address{Alfr\'ed R\'enyi Institute of Mathematics,
Hungarian Academy of Sciences,
Re\'altanoda utca 13-15, H-1053, Budapest, Hungary \newline
 \hspace*{4mm} ELTE - University of Budapest, Dept. of Geometry, Budapest, Hungary \newline \hspace*{4mm}
BCAM - Basque Center for Applied Math.,
Mazarredo, 14 E48009 Bilbao, Basque Country – Spain}
\email{nemethi.andras@renyi.mta.hu }
\thanks{D.K. was partially supported by the grant FP7-People-MCA-CIG, 334347.
Part of the work was done during D.K.'s posdoctoral fellowship in the University of Toronto.}
\thanks{A.N. was partially supported by NKFIH Grant  112735 and
ERC Adv. Grant LDTBud of A. Stipsicz at R\'enyi Institute of Math., Budapest}
\thanks{We thank R. Piene,  E. Shustin,  B. Sturmfels, J.J. Yang
for important advices.}


\begin{document}\setcounter{secnumdepth}{6} \setcounter{tocdepth}{1}

\begin{abstract}
Consider a space $X$ with the singular locus of positive dimension,
$Z=Sing(X)$. Suppose both $Z$ and $X$ are locally complete
intersections at each point.  The  transversal type of $X$ along $Z$
is generically constant but at some points of $Z$ it degenerates.
 In \cite{Kerner.Nemethi.Discriminant.Local}
 we have introduced (locally) the discriminant of the transversal type, a subscheme $\Db\sset Z$, that reflects these degenerations
  whenever the generic transversal type is ``ordinary". We have established the basic local properties of $\Db$.

In the current paper we consider the global case. We compute the
equivalence class of the discriminant in the Picard group, $[\Db]\in\Pic(Z)$.
If $X$ is a hypersurface, the discriminant is naturally stratified by the singularities of fibres in the projectivized normal cone $\P\cN_{X/Z}$.
 In this case (under some additional assumptions) we compute
 the classes of low codimension strata in the Chow group,
 $[\overline{\Si_{A_2}}],[\overline{\Si_{A_1A1}}]\in A^2(Z)$.

As immediate applications, we (re)derive
 the multi-degrees of the classical discriminant of projective complete intersections and
 bound the jumps of multiplicity of $X$ along $Z$ (when the singular locus is one-dimensional).
\end{abstract}
\maketitle\tableofcontents

\section{Introduction}

\subsection{The setup}\label{Sec.Introduction.Setup}
Let $\k$ be an algebraically closed field of zero characteristic, e.g. $\k=\C$. Let $X$ be either a reduced algebraic scheme or
 (for $\k=\C$)  a reduced complex-analytic space. We assume that the singular locus $Z:=Sing(X)$ is of positive dimension and connected.
 This second assumption is not restrictive; indeed, if $Sing(X)$ has several
connected component then we fix one of them, say $Z$,   and we replace $X$ by some neighborhood of $Z$.
 We always take $Z$ with its reduced  structure.

In many examples of non-isolated singularities one observes the following pattern.
 Fix a smooth point $\o\in Z$, and choose
 an embedding $(X,\o)\sset(M,\o)$ for some smooth germ  $(M,\o)$. Take
  a smooth germ, $(L^\bot,\o)\sset(M,\o)$,  transversal to $(Z,\o)$, such that $(L^\bot,\o)\cap(Z,\o)=\{\o\}$.
The singularity $(L^\bot\cap X,\o)$ is usually isolated and its type
 is in some sense generically constant along $Z$ (thus it is called the ``transversal singularity type").
 The points where the transversal singularity type degenerates usually form a subset of codimension 1 in $Z$,
 which can be defined locally.
  It is natural to call this subset the {\em discriminant of the transversal type}, $\Db\sset Z$.
Under certain assumptions (what we will also impose) regarding the pair $(Z,\o)\subset (X,\o)$, the space $\Db$ is even independent of the choice of $M$.

\subsection{The discriminant and its basic properties}\label{Sec.Introduction.Assumptions}
In \cite{Kerner.Nemethi.Discriminant.Local} we have rigorously defined  and studied this discriminant,  $\Db\sset Z$,
the subscheme which  reflects the degeneration of the transversal type. This was done under the following assumptions.
(See \S\ref{Sec.Background.Discriminant} for the precise definitions, examples and properties.)
\bee[i.]
\item The space $Z$ is a locally complete intersection at each point (l.c.i.).
\item For each point $\o\in Z$ the germ $(X,\o)$ is  a {\em strictly complete intersection} over $(Z,\o)$ (s.c.i.).
This is a slight strengthening of the notion of locally complete intersection, needed to ensure that the
the blowup $Bl_Z(X)$ is again a locally complete intersection at each point.
\item The transversal type of $X$ along $Z$  is  generically ``ordinary". Namely: there exists an open dense subset $\cU\sset Z$ such that for any
 point $\p\in\cU$, and any generic transversal section, $(L^\bot,\o)\cap (Z,\o)=\o$, the projectivized tangent cone, $\P T_{(L^\bot\cap X,\o)}$,
 is a smooth complete intersection of expected dimension.
\eee
Note that  being s.c.i. implies in particular the following fact as well: if $Z$  has several irreducible
components then the   multidegree  of $\P T_{(L^\bot\cap X,\o)}$ at
{\em generic} points of each  component is the same (see example \ref{Ex.SCI.hypersurface}).

We have established some of the basic properties of the subscheme $\Db\sset Z$, in particular:
\bei
\item $\Db$ is an effective  Cartier divisor in $Z$, fully determined locally and
 by the formal neighborhood of $Z$ in $X$;

\item
$\Db$ deforms flatly under those flat deformations of $X$ that preserve the multi-degree of $X$ along $Z$. More precisely:
 given a flat family, $\{X_t\}_{t\in(\k^1,0)}$, suppose that
 the singular locus is preserved, $Z=Sing(X_t)$.
Suppose each fibre $X_t$ satisfies the assumptions of \S\ref{Sec.Introduction.Assumptions}.
 Then the family $\{\Db(X_t)\}_{t\in(\k^1,0)}$ is flat and the
 class in the Neron-Severi group is preserved, $[\Db(X_t)]=const\in NS(Z)$.
\eei

\subsection{The equivalence class of $\Db$}\label{Sec.Introduction.Cohom.Class.Discrim}
In this paper we compute the class of the discriminant in the Picard group, $[\Db]\in Pic(Z)$.
In order to do this, we will embed $(X,\o)$ in some smooth space $M$, and the expression for
$[\Db]$ will use the embedding data as well. (This is the analogue of the expression of Euler characteristic of
a projective smooth complete intersection in terms of the degrees.)

More precisely, we will consider the following setup.
We fix an embedding $Z\sset X\sset M$, where the space $M$ is smooth at all the points of $X$.
  The normal sheaf, $\cN_{Z/M}$, of $Z$ in $M$ is locally free, we denote   its determinantal bundle by $det(\cN_{Z/M})$.
 Denote the codimension of $Z$ in $M$ by $k$.
   Suppose that $X$ is a globally complete intersection of codimension $r$, i.e. $X=\capl^r_{j=1}X_j\sset M$, for some hypersurfaces $X_j$.
Denote the generic multiplicity of $X_j$ along $Z$ by $p_j$.
 For each hypersurface $X_j\sset M$ take the class $[X_j]\in Pic(M)$  and restrict it to $Z$.
 We denote the restriction by the same letter, $[X_j]\in Pic(Z)$.

Furthermore,
we use the following combinatorial expression  related to the complete homogeneous symmetric polynomial
\beq\label{Eq.Combinatorial.Factor.Definition}
\cS_N:=\suml_{\substack{\{k_i\ge0\}\\\sum k_i=N}}\prodl^r_{i=1}(p_i-1)^{k_i}\stackrel{\S\ref{Sec.Background.Combinatorics}}{=\!=\!=}
\suml^r_{i=1}\frac{(p_i-1)^{N+r-1}}{\prodl_{j\neq i}(p_i-p_j)}.
\eeq
 (On the right hand side one could expect poles, when $p_i=p_j$, but these cancel, see \S\ref{Sec.Background.Combinatorics}.)

\bthe\label{Thm.Discriminant.Cohomology.Class}
Under the assumptions of \S\ref{Sec.Introduction.Assumptions}, the  class of $\Db\sset Z$,
  $[\Db]\in Pic(Z)$,
  is obtained as the coefficient of $t^{k-r}$ of the Taylor series at $t=0$ of the following expression:
\beq\label{Eq.Discriminant.Cohom.Class.First.Form}
[\Db]=
Coeff_{t^{k-r}}\Bigg[(1-t)^k\cdot\Big(\suml_j\frac{[X_j]/p_j}{1-p_jt}-\frac{[det(\cN_{Z/M})]}{1-t}\Big)
\cdot\prodl^r_{i=1}\frac{p_i}{1-p_it}\Bigg]\in Pic(Z),
\eeq
or explicitly:
\begin{multline}\label{Eq.Discriminant.Cohom.Class.Another.Form}
[\Db]=
(\prodl^r_{i=1} p_i)\Bigg(\suml_{j=1}^r\frac{[X_j]}{p_j}\suml^{k-r}_{l=0}(p_j-1)^l\cS_{k-r-l}-[det(\cN_{Z/M})]\cS_{k-r}\Bigg)=
\\=
(\prodl^r_{i=1}p_i)\Bigg(\suml^r_{j=1}\frac{[X_j]}{p_j}\suml^r_{i=1}\frac{(p_j-1)^k-(p_i-1)^k}{(p_j-p_i)\prodl_{l\neq i}(p_i-p_l)}
-[det(\cN_{Z/M})]\suml^r_{j=1}\frac{(p_j-1)^{k-1}}{\prodl_{i\neq j}(p_j-p_i)}\Bigg).
\end{multline}
\ethe

\beRs\label{Remarks}
\bee[i.]

\item In this theorem we assume $X\sset M$ a {\it globally}
 complete intersection. The generalization to the case of {\it locally} complete intersection
  is done in  \S\ref{Sec.X.locally.Complete.Intersection}. The next discussion aims to motivate the form
  of this generalization. First we   order the sequence of multiplicities as
\[
p_1=\cdots=p_{k_1}<p_{k_1+1}=\cdots=p_{k_2}<p_{k_2+1}=\cdots p_r.
\]
Then  we rewrite the sum $\suml^r_{j=1}\frac{[X_j]}{p_j}\suml^{k-r}_{i=1}(\dots)$ in the last equation
    (\ref{Eq.Discriminant.Cohom.Class.Another.Form})
in the form
$\suml_l\frac{[X_{k_l+1}]+\cdots+[X_{k_{l+1}}]}{p_{k_l}}\suml^{k-r}_{i=1}(\dots)$.
This suggests that even if $X\subset M$ is not a globally complete intersection, there exists
a filtration  of $\cN_{X/M}$ by a sequence of subbundles,
 $\cN_{X/M}=\cN_{p_{k_0}}\supseteq \cN_{p_{k_1}}\supseteq\cdots\supseteq\cN_{p_{k_r'}}=\{0\}$, such that
 $\quotients{\cN_{p_{k_l}}}{ \cN_{p_{k_{l+1}}}}$ has rank $k_{l+1}-k_l$, and
  whenever $X\subset M$ is globally complete intersection then
 $[det(\quotients{\cN_{p_{k_l}}}{ \cN_{p_{k_{l+1}}}})]=[X_{k_l+1}]+\cdots+[X_{k_{l+1}}]$.
  In \S\ref{Sec.X.locally.Complete.Intersection} we construct this filtration when $X\sset M$ is a locally complete intersection. In fact, it is convenient to index the filtration of subbundles as
  $\cN_{X/M}=\cN_{p_0}\supseteq \cN_{p_1}\supseteq\cdots\supseteq\cN_{p_r}=\{0\}$, such that the quotients
   $\quotients{\cN_{p_{k_l}}}{\cN_{p_{k_l+1}}}$ has rank  $k_{l+1}-k_l$, while the other
   quotients are trivial.
   Then one obtains
\beq\label{Eq.Discriminant.Cohom.Class.loc.complete.inters}
[\Db]=
(\prod p_i)\Bigg(\suml_{j=1}^r\frac{[det(\quotients{\cN_{p_{j-1}}}{\cN_{p_{j}}})]}{p_j}\suml^{k-r}_{l=0}(p_j-1)^l\cS_{k-r-l}-[det(\cN_{Z/M})]\cS_{k-r}\Bigg).
\eeq
\item
As we see from the formulas, ``in most cases" $[\Db]\neq0$, i.e. the transversal type necessarily degenerates in codimension one.
 Thus theorem \ref{Thm.Discriminant.Cohomology.Class} is an {\em obstruction} to the
 naive expectation (that might arise from the real differential geometry): ``in the (very) generic case the transversal type of $X$ along $Z$ is generic  and does not degenerate".
  This obstruction is a far-reaching extension of Proposition on page 3 of \cite{Aluffi-1996}.
\item
Vice versa, as  $\Db\sset Z$ is an effective divisor,  $[\Db]=0$ means $\Db=\empty$. Therefore, if the data $\{p_j\}$, $[X_j]$, $[det(\cN_{Z/M})]$
 has been ``fine tuned" to ensure the vanishing $[\Db]=0$, then
 the transversal type of $X$ is ordinary {\it everywhere} along $Z$
 for {\em any} germ $(X,Z)$ that realizes this data and has (a priori)
  {\it generically} ordinary transversal type.
\item As is mentioned above, the subscheme $\Db\sset Z$ is fully determined by the formal neighborhood of $Z$ in $X$, \cite[\S 5.2]{Kerner.Nemethi.Discriminant.Local}.
In the above formulae, we express $[\Db]$ in terms of the embedding data.  By
 Theorem \ref{Thm.Discriminant.Cohomology.Class},  the class $[\Db]\in \Pic(Z)$ is determined by the normal sheaves, $\cN_{Z/M}$, $\cN_{X/M}$,
  i.e., by the first infinitesimal neighborhoods of $Z,X$ in $M$.
  (However, by the nature of the local construction of the discriminant, it is independent of the choice of $M$.)
  \item
When $X\sset M$ is a hypersurface, the formula simplifies into
\beq\label{Eq.Discrim.Class.for.Hypersurfaces}
[\Db]=(p-1)^{k-1}\Big(k[X]- p\cdot [\det(\cN_{Z/M})]\Big).
\eeq
By the adjunction formula, $w_Z=w_M\otimes\cO_Z\otimes det(\cN_{Z/M}^\vee)$, this can be rewritten in the form:
\beq\label{eq:deJong}
[\Db]=(p-1)^{k-1}\Big(k[X]+p(K_M|_Z-\omega_Z)\Big),
\eeq
where $K_M$ and $\omega_Z$ are the dualizing sheaves.
Here, the right hand side is meaningful even for an arbitrary reduced Gorenstein scheme $Z$. Thus the natural task is to extend the whole
 setting to this generality. This was done in \cite{de Jong-de Jong} for the particular case, when $X\sset M$ is a hypersurface in a complex manifold,
  $p=2$, and $Z$ is a  Gorenstein curve. If moreover $Z$ is locally complete intersection at each point then
  our last equation (\ref{eq:deJong}) reproves the main
   theorem of \cite{de Jong-de Jong}.
\eee
\eeRs

\subsection{The proof Theorem \ref{Thm.Discriminant.Cohomology.Class} and the
 relation with the multi-weighted degrees of the classical discriminant}\label{Sec.Introduction.Proof.and.Further.Results}
It turns out that  the statement of Theorem \ref{Thm.Discriminant.Cohomology.Class}
implies the formula for the  multi-weighted degrees of the discriminant
 of projective complete intersections and, at the same time, it can be reproved based on this formula. This fact is reflected totally in the proof, which runs
  in three steps, see \S\ref{Sec.Proofs}.
\bee[i.]
\item When $Z$ is algebraic and smooth, and the blowup $Bl_Z(X)$ is smooth too,
  we compute the equivalence class in the Chow group, $[\Db]\in A^1(Z)$, using Thom-Porteus formula, see
   \S\ref{Sec.Proofs.Cohomology.Class.Discriminant.Thom.Porteus}.

\item Then, using the class $[\Db]\in A^1(Z)$, we compute the multi/weighted-degrees of the classical discriminant of complete intersections,
 see \S\ref{Sec.Proofs.Classical.Discriminant.Multi-Degree}
 In the hypersurface case (the discriminant of one homogeneous form) these degrees are well known, see e.g.  \cite{Gelfand-Kapranov-Zelevinsky}.
  For complete intersections, it seems, they have been computed for the first time in \cite{Benoist}. Our derivation, using $[\Db]\in A^1(Z)$, is
   independent of these works.
\item Finally, using the multi/weighted-degrees of the classical discriminant, we compute the class $[\Db]\in Pic(Z)$, or, rather, we identify the
 line bundle $\cO_Z(\Db)$, in the required general case of $X$, see \S\ref{Sec.Proofs.Cohomology.Class.Discriminant.Classical.Style}.
 In particular,
  if one just uses/accepts the (known) multi/weighted-degrees then this last step
  provides a simple and short  proof of Theorem \ref{Thm.Discriminant.Cohomology.Class},
  totally independent of the Thom-Porteus formula.

\eee

\subsection{The equivalence classes of the strata  $\lSi_{A_2}$, $\lSi_{A_1,A_1}$}\label{Sec.Introduction.Equivalence.Classes.Strata}
Any stratification of the classical discriminant (in the parameter space of projective hypersurfaces or projective complete intersections)
 induces a stratification of $\Db$, see \cite[\S 4.3]{Kerner.Nemethi.Discriminant.Local}.  In many cases the generic points of $\Db$ correspond
  to a nodal, $A_1$--degeneration of the fibres in projectivized normal cone $\P\cN_{X/Z}$.
 Similarly, the codimension-two strata (in $Z$) correspond to a cusp, $A_2$, or to two nodes, $A_1A_1$.
In \S\ref{Sec.Cohomol.Classes.of.Further.Strata}, we compute the classes of these
codimension-two strata $[\lSi_{A_2}]$, $[\lSi_{A_1,A_1}]$ in the Chow group $A^2(Z)$,
 provided $X\sset M$ is a hypersurface and $Z$ is smooth. This computation uses the Thom  polynomial technique of \cite{Kazarian2000}, \cite{Kazarian2006}.

\subsection{}
In \S\ref{Sec.Application.Bound.on.jump.of.multiplicity}, we consider the case of the one-dimensional
singular locus. In this case the discriminant is a collection of points, thus we can
bound the total
``amount of degenerations" along the singular locus. We give an explicit upper bound on the total jumps
of multiplicity of $X$ along $Z$.

\subsection{History}\label{Sec.Introduction.History}
\bee[i.]
\item The discriminant of transversal singularity type appears naturally
in geometry and singularity theory. In some particular cases
 this was considered already by Salmon, Cayley, Noether and Zeuthen, see \cite{Piene1978}.

 The class of $\Db$ for projective surface,
$X\sset\P^3$, with ordinary singularities along the singular locus, a curve, goes back (probably) to
the early history. For a computation of the degree of $\Db$ and various other invariants see \cite{Piene1978} and \cite[example 9.3.7]{Fulton-book}

\item The case of one-dimensional singular locus, i.e. when $Z$ is a curve, with generic transversal type $A_1$,
 was thoroughly studied by Siersma, see e.g. \cite{Siersma2000}. The local degree of the discriminant, called also
 ``the virtual number of $D_\infty$ points" was studied in \cite{Pellikaan-PhD}, \cite{Pellikaan} and \cite{de Jong}.
 In particular, they show the pathological behavior of $\Db$ when $Z$ is
 not a locally complete intersection. In \cite{de Jong-de Jong} the degree of $[\Db]$ is computed for the
 case when $X\sset M$ is a projective hypersurface, $Z=Sing(X)$ is of (pure) dimension one and the generic transversal type is $A_1$.
\eee

\section{Preparations}
\subsection{Some  combinatorics related to $S_N$}\label{Sec.Background.Combinatorics}
The formulae are expressed in terms of the complete homogeneous symmetric polynomial
 $\cS_N(e_1,\dots,e_r):=\suml_{\substack{\{k_i\ge0\}\\\sum k_i=N}}\prodl^r_{i=1}e_i^{k_i}$.
In our presentation we use the following properties of this polynomial.
(Some proofs are indicated, we invite the reader to check the validity
of the other identities as well.)
\bee[i.]
\item
$\cS_N(e_1,\dots,e_r)=\suml^r_{i=1}\frac{e^{N+r-1}_i}{\prodl_{j\neq i}(e_i-e_j)}$, see e.g. \cite[ex.7.4, page 450]{Stanley}.
\item
For $q<r-1$ one has $\suml^r_{i=1}\frac{e^q_i}{\prodl_{j\neq i}(e_i-e_j)}=0$. This follows, e.g. by considering Lagrange interpolation polynomial,
  $L_{\{y_i\}}(x)=\suml^r_{i=1}\prodl^r_{\substack{j=1\\j\neq i}}\frac{x-x_j}{x_i-x_j}y_i$. It is of degree $<r$ and satisfies
   $\{L_{\{x^{q}_i\}}(x_j)=x^{q}_j\}_{j=1,\dots, r}$. Thus one has $L_{\{x^{q}_i\}}(x)=x^{q}$. Finally, extract the coefficient of $x^{r-1}$ in
 $L_{\{x^{q}_i\}}(x_j)$.

\item $\suml^r_{i=1}\frac{e^q_i-e^q_j}{e_i-e_j}\frac{1}{\prodl_{l\neq i}(e_i-e_l)}=0$ for any fixed $j$ and $0\le q<r$.
 Here expand $\frac{e^q_i-e^q_j}{e_i-e_j}=\sum_l e^l_ie^{q-l-1}_j$ and use part ii.
\item $\suml^{k-r}_{q=0}e_j^qS_{k-r-q}=\suml^r_{i=1}\frac{e^{k}_j-e^{k}_i}{e_j-e_i}\frac{1}{\prodl_{l\neq i}(e_i-e_l)}$ \ for any $k\ge r$. 
\item  $\suml^r_{\substack{i,j=1\\i\neq j}}\frac{e^k_i-e^k_j}{e_i-e_j}\frac{1}{\prodl_{l\neq i}(e_i-e_l)}=0$ \ for any $k\ge0$. 
\item  $\suml^r_{i=1}\Big(\suml^r_{j=1}\frac{e^k_i-e^k_j}{e_i-e_j}\frac{1}{\prodl_{l\neq j}(e_j-e_l)}\Big)=
\suml^r_{i=1}\frac{ke^{k-1}_i}{\prodl_{l\neq i}(e_i-e_l)}$ \ for any $k\ge0$. 
\eee

\subsection{Some relevant facts about the discriminant of transversal type}\label{Sec.Background.Discriminant}
We recall briefly the definition and properties of the subscheme $\Db\sset Z$, see \S4 and \S5 of \cite{Kerner.Nemethi.Discriminant.Local}.

\subsubsection{The normal cone}\label{Sec.Background.Normal.Cone}
A  germ $(Z,\o)\sset(\k^N,\o)$ defines the ``associated graded" functor that acts on the ring, ideals and modules,
 e.g.
\beq
 \cO_{(\k^N,\o)}\stackrel{gr_{(Z,\o)}}{\rightsquigarrow}gr_{(Z,\o)}\cO_{(\k^N,\o)}=\oplusl_{j\ge0}\quotient{I_{(Z,\o)}^j}{I_{(Z,\o)}^{j+1}}.
\eeq
In our case $(Z,\o)$ is a reduced complete intersection. Fix a
  regular sequence of generators, $I_{(Z,\o)}=(g_1,\dots,g_k)\sset\cO_{(\k^N,\o)}$. This defines the isomorphism of algebras,
   $gr_{(Z,\o)}\cO_{(\k^N,\o)}\approx \cO_{(Z,\o)}[\uy]$, here $\uy=(y_1,\dots,y_k)$.

To any element $f\in \cO_{(\k^N,\o)}$ we associate its leading term, $\tf\in gr_{(Z,\o)}\cO_{(\k^N,\o)}$, as follows. Suppose $f$ is of
 order $p$ on $(Z,\o)$, i.e. $f\in I_{(Z,\o)}^p\smin I_{(Z,\o)}^{p+1}$. Accordingly, we expand  $f$ in powers of generators of $I_{(Z,\o)}$,
 \beq\label{Eq.Expansion.of.f.in.g}
f=\suml_{\sum m_j=p}g^{m_1}_1\dots g^{m_k}_k a_{m_1\dots m_k}, \ \text{ for some elements } \ \{a_{m_1\dots m_k}\in  \cO_{(\k^N,\o)}\}.
\eeq
Then the leading term of $f$ is: $\tf=\sum y^{m_I}_I a_{m_I}|_{(Z,\o)}\in\cO_{(Z,\o)}[\uy]$.

 The coefficients $\{a_{m_I}\}$ are not unique, e.g. one has the Koszul relations among $\{g_i\}$.
Nevertheless, as $(Z,\o)$ is a complete intersection,  the restrictions $\{a_{m_I}|_{(Z,\o)}\in \cO_{(Z,o)}\}$ are defined uniquely,
 see \cite[\S2.4.1]{Kerner.Nemethi.Discriminant.Local}.

\

To the germ $(X,\o)\sset(\k^N,\o)$ we associate the ideal $gr_{(Z,\o)}I_{(X,\o)}\sset gr_{(Z,\o)}\cO_{(\k^N,\o)}$, generated by all the leading
 terms of elements of $I_{(X,\o)}$. This $gr_{(Z,\o)}I_{(X,\o)}$ is the conormal ideal of $(Z,\o)$ in $(X,\o)$. Accordingly, the normal cone of
  $(Z,\o)$ in $(X,\o)$ is defined as
\beq\label{Eq.Expansion.to.define.normal.cone}
  \cN_{(X,\o)/(Z,\o)}=V(gr_{(Z,\o)}I_{(X,\o)})\sset Spec(gr_{(Z,\o)}\cO_{(\k^N,\o)})=(Z,\o)\times Spec(\k[\uy]).
\eeq

\subsubsection{The strictly complete intersections}\label{Sec.Background.s.c.i}
The germ $(X,\o)$ is called a \underline{strictly complete intersection} over $(Z,\o)$, abbreviated s.c.i., if the normal cone $  \cN_{(X,\o)/(Z,\o)}$ is a
 complete intersection over $(Z,\o)$. This is a strengthening of the ordinary property of complete intersection,
  see \cite[\S2.5]{Kerner.Nemethi.Discriminant.Local}.

\bex (See Example 2.10 of \cite{Kerner.Nemethi.Discriminant.Local})
Suppose  $(X,\o)$ is a hypersurface singularity, with $(Z,\o)=Sing(X,\o)$.
 If $(Z,\o)$ is irreducible  then $(X,\o)$ is s.c.i. over $(Z,\o)$. If $(Z,\o)$
 is reducible then  $(X,\o)$ is s.c.i. over $(Z,\o)$ iff the generic multiplicity of $(X,\o)$ is the same on all the components of $(Z,\o)$.
\eex

The geometric interpretation of being s.c.i. is the following: $(X,\o)$ is a complete intersection and for the blowup,
  as on the diagram \eqref{Eq.Diagram.for.Critical.Locus},
 the subscheme $\tX\cap E\sset E$  is a complete intersection over $(Z,\o)$.

\

If $(X,\o)$ is s.c.i. over $(Z,\o)$, then the ideal $I_{(X,\o)}$ has a basis of special type, called \underline{a good basis}, $I_{(X,\o)}=(f_1,\dots,f_r)$,
 such that
 the leading terms of $\{f_i\}$ generate the associated graded ideal,
 $gr_{(Z,\o)}I_{(X,\o)}=(\tf_1,\dots,\tf_r)$, \cite[\S2.6]{Kerner.Nemethi.Discriminant.Local}.
The \underline{multiplicity sequence} of $(X,\o)$ along $(Z,\o)$ is the set of multiplicities of $f_i$ along $(Z,\o)$. Though the good basis is not unique,
 the multiplicity sequence does not depend on  the choice of a good basis, see Proposition 2.19 of \cite{Kerner.Nemethi.Discriminant.Local}.

\subsubsection{Generically ordinary singularities}
Set $Z=Sing(X)$ and let $\o\in Z$ be a point. Take a smooth transversal germ of complementary dimension, $(L^\bot,\o)\cap(Z,\o)=\o$.
 (The transversality here means that
  the tangent cones $T_{(L^\bot,\o)}$, $T_{(Z,\o)}$ are mutually generic.)
 The $\mu=const$ type of this intersection depends on $\o$, sometimes it depends even on the choice of $L^\bot$, see
  Example 1.1 in \cite{Kerner.Nemethi.Discriminant.Local}, however, for generic $\o$ the following definition is meaningful.

 $X$ is called \underline{generically ordinary} along $Z$ if for generic $\o$ and generic $L^\bot$ the projectivized tangent cone, $\P T_{(X\cap L^\bot,\o)}$,
  is a smooth complete intersection of expected dimension. This condition can be formulated algebraically, in terms of the localizations
   of the projectivized normal cone, $\P\cN_{(X,\o)/(Z,\o)}$.
\bex\label{Ex.SCI.hypersurface}
Let $X\sset M$ be a hypersurface and $Z=Sing(X)$. Then $X$ is generically ordinary along $Z$ if for the generic point $\o$ on any irreducible component of $Z$
 and the generic transversal slice $(L^\bot,\o)$, the hyperurface singularity $(X\cap L^\bot,\o)\sset(L^\bot,\o)$ is an ordinary multiple point.
  Namely, it can be defined by equation $\{f_p(\ux)+f_{>p}(\ux)=0\}$, where the homogeneous polynomial $f_p(\ux)$ is non-degenerate.
\eex

\subsubsection{The definition of $\Db$ as the pullback of the classical discriminant}
Equation \eqref{Eq.Expansion.of.f.in.g} associates to a good basis $I_{(X,\o)}=(f_1,\dots,f_r)\sset\cO_{(\k^N,o)}$ the collection of coefficients,
  $\{a^{(j)}_{m_I}|_{(Z,\o)}\}$.
 This collection of coefficients
 defines a subscheme in $\P^{k-1}\times (Z,\o)$. If $X$ is generically ordinary along $Z$, then, for generic $\o\in Z$, this subscheme will
  be a family, over $(Z,\o)$, of smooth   projective complete intersections.
 Thus we get a (rational) map, from  $Z$ to the parameter space of projective complete intersections.
The image of this map does not lie fully inside the classical discriminant $\De\sset\prodl^r_{j=1} |\cO_{\P^k-1}(p_j)|$. Using this map we pullback the
 classical discriminant of complete intersections. By definition, $\Db$ is this pullback, see \S4.1 of \cite{Kerner.Nemethi.Discriminant.Local}.

\subsubsection{The definition of $\Db$ as the pushforward of the critical locus}

An equivalent definition is as the image of the critical locus, \cite[\S4.2]{Kerner.Nemethi.Discriminant.Local}, as follows. In the embedded case, $Z\sset X\sset M$, blowup over $Z$, as on the
 diagram (here $E$ is the exceptional divisor and $\tX\sset Bl_Z(M)$ is the strict transform of $X$):
\beq\label{Eq.Diagram.for.Critical.Locus}
\bM Crit(\pi)\sset \P\cN_{X/Z}=
 \tX\cap E\stackrel{\nu}{\hookrightarrow}&E&\approx \P\cN_{Z/M}&\sset Bl_Z(M)
 \\
\hspace{3.5cm}\pi\searrow
&\downarrow&&\downarrow
\\\hspace{3cm}&Z&\sset& M
\eM
\eeq
(Note that both $Z$ and $\tX\cap E$ are in general non-smooth, but all the singularities are at worst locally complete intersections.)

The critical locus of the projection $\pi$ is defined, as usual, by the Fitting ideal of
 the relative cotangent sheaf,
 \beq
 I(Crit(\pi))=Fitt_{k-r}(\Om^1_{\tX\cap E/Z}),\quad \text{here}\quad 0\to \pi^*\Om^1_Z\to \Om^1_{\tX\cap E}\to \Om^1_{\tX\cap E/Z}\to0.
\eeq
  Then, at the points where the restriction $\pi|_{Crit(\pi)}$
  is finite, the discriminant $\Db$ is defined as the direct image of the critical locus, with the Fitting scheme structure: $I_{\Db}=Fitt_0(\pi_*\cO_{(Crit(\pi))})$.
One then takes the topological closure to get the full discriminant, $\Db\sset Z$.

\section{The proof of theorem \ref{Thm.Discriminant.Cohomology.Class}}\label{Sec.Proofs}

\subsection{Computation of the class $[\Db]\in A^1(Z)$ using Thom-Porteus formula}\label{Sec.Proofs.Cohomology.Class.Discriminant.Thom.Porteus}
\bthe\label{Thm.Discriminant.Cohomology.Class.via.Thom.Porteus}
We work under the assumptions of \S\ref{Sec.Introduction.Assumptions}. Suppose $X$ is a globally complete intersection, $X=\capl^r_{j=1}X_j\sset M$.
Assume moreover that $Z,\tX\cap E$ are smooth algebraic schemes.
 Then the class $[\Db]\in A^1(Z)$  is given by equations \eqref{Eq.Discriminant.Cohom.Class.First.Form} and \eqref{Eq.Discriminant.Cohom.Class.Another.Form},
  where now $[X_j]\in A^1(Z)$ and $[det(\cN_{Z/M})]$ is replaced by the first Chern class $c_1(\cN_{Z/M})$.
\ethe
\bpr
{\bf Step 1.} We work in the setup of the diagram in equation \eqref{Eq.Diagram.for.Critical.Locus}.
Let $I=I_{Z/M}\sset \cO_M$ be the ideal sheaf of $Z$ in $M$.

The fibers of the projection $\tX\cap E\proj Z$ are
 projective varieties, generically (over $Z$) they are complete intersections of dimension $(k-r-1)$.

The discriminant is the pushforward of the critical locus, $\Db=\pi_*Crit(\pi)$, with the Fitting scheme structure,
 i.e. $\cO_{\Db}=\pi_*\cO_{Crit(\pi)}\sset\cO_Z$ and
 $I_{\Db/Z}=Fitt_0(\pi_*\cO_{Crit(\pi)})$. Thus the equivalence class of $\Db$ is obtained by the pushforward, $[\Db]=\pi_*[Crit(\pi)]\in A^1(Z)$.

The subscheme $Crit(\pi)\sset\tX\cap E$ is defined as the degeneracy locus of the map $T_{\tX\cap E}\to \pi^* T_Z$. The  scheme
 structure is given  by the Fitting ideal,  $Fitt_{k-r}(\Om^1_{\tX\cap E/Z})$, \cite[\S3.3]{Kerner.Nemethi.Discriminant.Local}.
  Then the class $[Crit]\in  A^{k-r}(\tX\cap E)$ is obtained by Thom-Porteus formula, \cite[\S14.4]{Fulton-book}:
\beq
[Crit(\pi))]=(-1)^{k-r}c_{k-r}(T_{\tX\cap E}-\pi^*T_Z).
\eeq
Here the total Chern polynomials of a vector bundle is $c_t(V)=\sum_{i\ge0}t^ic_i(V)$, while
the relative Chern class is determined by the expansion,
$c_t(T_{\tX\cap E}-\pi^*T_Z)=\frac{c_t(T_{\tX\cap E})}{c_t(\pi^*T_Z)}$.

For this we use the adjunction sequence $0\to T_{\tX\cap E}\to T_E|_{\tX\cap E}\to\cN_{\quotients{\tX\cap E}{E}}\to0$, where the later bundle is
the normal bundle of $\tX\cap E$ in $E$.
Recall that $\tX\cap E\sset E$ is a complete intersection, thus $c_t(\cN_{\quotients{\tX\cap E}{E}})=\prod(1+t[\tX_j])$,
where on $E$ we have: $[\tX_j]=[\pi^*X_j]+p_jc_1(\cO_E(1))$, for the tautological bundle $\cO_E(1)$.
 (Here $[\pi^*X_j]$ is the restriction of the class $[\pi^*X_j]\in A^1(Bl_Z(M))$ to $E$.)

The Chern class of the relative tangent bundle $T_{E/Z}=T_E-\pi^*T_Z$ can be computed as follows
An element of $T_{E/Z}$ at each point over $\o\in Z$ corresponds to a two dimensional vector subspace
in $\cN_{Z/M}|_{\o}$. Such elements are produced as maps from the tautological line bundle on $E$ to $\pi^*\cN_{Z/M}$,
 so $Hom(\cO_E(-1),\pi^*\cN_{Z/M})\to T_{E/Z}\to0$. The kernel of this surjection consists of maps that send
 the line to itself, i.e.
\beq
 0\to Hom(\cO_E(-1),\cO_E(-1))\to Hom(\cO_E(-1),\pi^*\cN_{Z/M})\to T_{E/Z}\to0.
\eeq
This gives the exact sequence, see also \cite[Example 3.2.11]{Fulton-book}:
\beq
0\to\cO_E\to\pi^*(\cN_{Z/M})\otimes\cO_E(1)\to T_E\to\pi^*T_Z\to0.
\eeq
 Thus $\frac{c_t(T_E)}{c_t(\pi^*T_Z)}=c_t\Big(\pi^*(\cN_{Z/M})\otimes\cO_E(1)\Big)=
\Big(1+tc_1(\cO_E(1))\Big)^k c_\tau\Big(\pi^*(\cN_{Z/M})\Big)$,
where $\tau=\frac{t}{1+tc_1(\cO_E(1))}$.
Altogether:
\beq\label{Eq.Critical.Locus.Cohom.Class}
[Crit(\pi)]=(-1)^{k-r}Coeff_{t^{k-r}}\Bigg(
\frac{\suml_{i=0}^k\Big(1+tc_1(\cO_E(1))\Big)^{k-i}t^jc_i(\pi^*(\cN_{Z/M}))}{\prodl^r_{j=1}(1+t[\tX_j])}\Bigg)\in
A^{k-r}(\tX\cap E).
\eeq

{\bf Step 2.}
The class $[\Db]\in A^1(Z)$ is obtained by pushforward $\pi_*[Crit(\pi)]$. To compute this we first obtain
the class $j_*[Crit(\pi)]\in A^k(E)$, see the diagram \eqref{Eq.Diagram.for.Critical.Locus}.
 For this we multiply the expression in equation \eqref{Eq.Critical.Locus.Cohom.Class} by $[\tX\cap E]=\prod[\tX_i]$.
Then we apply the Gysin homomorphism, i.e. the projection $A^k(E)\to A^1(Z)$.

Here is the computation.  For brevity we denote $c_1(\cO_E(1))$ by $-[E]$, which means: the class $-[E]\in Pic(Bl_Z(M))$ restricted to $E$. We have the class
\beq
A^{k}(E)\ni \nu_*[Crit(\pi)]=(-1)^{k-r}\prodl^r_{j=1}\Big([\pi^*X_j]-p_j[E]\Big)Coeff_{t^{k-r}}\Bigg(
\frac{\suml_{i=0}^k\Big(1-t[E]\Big)^{k-i}t^ic_i(\pi^*(\cN_{Z/M}))}{\prodl^r_{j=1}\Big(1+t[\pi^* X_j]-tp_j[E]\Big)}\Bigg).
\eeq
As the resulting class is in $A^1(Z)$ we expand the expression in $[\pi^*X_j]$ and $c_i$ and
keep only the terms of order $\le1$. In particular we discard $c_{i>1}$.
 Thus:
\beq
(-1)^{k-r}\prodl^r_{j=1}\Big([\pi^*X_j]-p_j[E]\Big)
\rightsquigarrow
(-1)^k(\prodl^r_{j=1} p_j)\Big([E]^r-[E]^{r-1}\suml^r_{j=1}\frac{[\pi^* X_j]}{p_j}\Big)
\eeq
and
\beq
\frac{\suml_{i=0}^k\Big(1-t[E]\Big)^{k-i}t^ic_i(\pi^*(\cN_{Z/M}))}{\prodl^r_{ji=1}\Big(1+t[\pi^* X_j]-tp_j[E]\Big)}
\rightsquigarrow
\frac{(1-t[E])^k\Big(1-\sum\frac{t[\pi^*X_j]}{1-p_jt[E]}\Big)+(1-t[E])^{k-1}tc_1(\pi^*(\cN_{Z/M})}{\prod(1-p_jt[E])}.
\eeq
Denote
\begin{multline}
B:=(-1)^k(\prodl^r_{j=1} p_j)\Big([E]^r-[E]^{r-1}\suml^r_{j=1}\frac{[\pi^* X_j]}{p_j}\Big)\cdot
\\
\cdot Coeff_{t^{k-r}}\Bigg(\frac{(1-t[E])^k\Big(1-\sum\frac{t[\pi^*X_j]}{1-p_jt[E]}\Big)+(1-t[E])^{k-1}tc_1(\pi^*(\cN_{Z/M})}{\prod(1-p_jt[E])}\Bigg).
\end{multline}
Applying the Gysin homomorphism amounts to the substitution
(after expanding in powers of $E$): $E^i\to(-1)^is_{i-k+1}$, where $s_{i-k+1}$ is the Segre class of $\cN_{Z/M}$.
Again, as the resulting class is in $A^1(Z)$ we need only the terms linear in $[X_j]$ or $s_1(\cN_{Z/M})$, but
no higher terms. Hence we only need to extract the coefficients of $E^k$ and $E^{k-1}$:
\beq\ber
(-1)^k s_1(\cN_{Z/M})Coeff_{[E]^k}(B)+(-1)^{k-1}Coeff_{[E]^{k-1}}(B)=(-1)^{k-1}Coeff_{[E]^k}\Big(c_1(\cN_{Z/M})B+[E]B\Big)=\\
=-(\prodl^r_{j=1}p_j)Coeff_{[E]^kt^{k-r}}\Bigg(c_1(\cN_{Z/M})\frac{[E]^r(1-[E]t)^k}{\prod(1-p_jt[E])}+
[E]^{r}(1-[E]t)^k\frac{E-(\sum\frac{[E]t[X_j]}{1-p_jt[E]})+\frac{[E]tc_1(\cN_{Z/M})}{(1-t[E])}
-\suml^r_{j=1}\frac{[X_j]}{p_j}}{\prod(1-p_jt[E])}\Bigg).
\eer\eeq
We note that $Coeff_{[E]^kt^{k-r}}\Big(\frac{[E]^{r+1}(1-[E]^kt)}{\prod(1-p_it[E])}\Big)=0$ so we omit this term.
The rest equals:
\beq\ber
-(\prodl^r_{j=1}p_j)Coeff_{[E]^{k-r}t^{k-r}}\Bigg(c_1(\cN_{Z/M})\frac{(1-[E]t)^k}{\prod(1-p_jt[E])}+
(1-[E]t)^k\frac{-(\sum\frac{[E]t[X_j]}{1-p_jt[E]})+\frac{[E]tc_1(\cN_{Z/M})}{(1-t[E])}
-\suml^r_{i=1}\frac{[X_i]}{p_i}}{\prod(1-p_jt[E])}\Bigg)
\\
\stackrel{\tau:=[E]t}{=\!=\!=}(\prodl^r_{j=1}p_j)Coeff_{\tau^{k-r}}\Bigg(-c_1(\cN_{Z/M})\frac{(1-\tau)^k}{\prod(1-p_j\tau)}-
(1-\tau)^k\frac{(-\sum\frac{\tau[X_j]}{1-p_j\tau})+\frac{\tau c_1(\cN_{Z/M})}{(1-\tau)}
-\suml^r_{j=1}\frac{[X_j]}{p_j}}{\prod(1-p_j\tau)}\Bigg)=\\=
Coeff_{\tau^{k-r}}\Bigg(\Big((1-\tau)^k\sum\frac{[X_j]/p_j}{1-p_j\tau}-(1-\tau)^{k-1}c_1(\cN_{Z/M})\Big)(\prodl^r_{j=1}\frac{p_j}{1-p_j\tau})
\Bigg)=\\=
Coeff_{\tau^{k-r}}\Bigg[(1-\tau)^k\Big(\sum\frac{[X_j]/p_j}{1-p_j\tau}-\frac{c_1(\cN_{Z/M})}{1-\tau}\Big)(\prodl^r_{j=1}\frac{p_j}{1-p_j\tau})
\Bigg].
\eer\eeq
This proves the formula of equation \eqref{Eq.Discriminant.Cohom.Class.First.Form}.

{\bf Step 3.}
Now we obtain the explicit form of the coefficients. The coefficient of $c_1(\cN_{Z/M})\prodl^r_{i=1}p_i$ is:
\begin{multline}
Coeff_{\tau^{k-r}}\frac{(1-\tau)^{k-1}}{\prodl_j(1-p_j\tau)}\stackrel{\tau\to -t}{=}(-1)^{k-r}Coeff_{t^{k-r}}\frac{(1+t)^{k-1}}{\prodl_j(1+p_jt)}=
(-1)^{k-r}Coeff_{t^{k-r}}\Big(\frac{(1+t)^{k-r-1}}{\prodl_j(1-\frac{(1-p_j)t}{1+t})}\Big)=
\\=
(-1)^{k-r}Coeff_{t^{k-r}}\Big(\sum_{\substack{\{k_j\ge0\}_j}}t^{\sum k_j}(1+t)^{k-r-1-\sum k_j}\prodl_j(1-p_j)^{k_j}\Big).
\end{multline}
 In the obtained expression,
if $\sum k_j\le k-r-1$ then there is no term of order $(k-r)$. If $\sum k_j\ge k-r$ then we have
$\frac{t^{\sum k_j}}{(1+t)^{\sum k_j-k+r+1}}=t^{\sum k_j}(1+\cdots)$. So, if $\sum k_j>k-r$ then all the terms here
are of order bigger than $(k-r)$. While for $\sum k_j=k-r$ there is precisely one term, with coefficient 1.
Therefore
\beq
(-1)^{k-r}Coeff_{t^{k-r}}\Big(\sum_{\substack{\{k_j\ge0\}_j}}t^{\sum k_j}
(1+t)^{k-r-1-\sum k_j}\prodl_j(1-p_j)^{k_j}\Big)=(-1)^{k-r}\sum_{\substack{k_j\ge0, \forall\ j\\\sum k_j=k-r}}
\prodl_j(1-p_j)^{k_j}=S_{k-r}.
\eeq
By similar arguments we get:
$Coeff_{\tau^{k-r}}\Big(\frac{(1-\tau)^k[X_j]/p_j}{(1-p_j\tau)}\prodl^r_{i=1}\frac{p_i}{1-p_i\tau}\Big)=
\frac{[X_j]}{p_j}(\prod p_i)\suml^{k-r}_{l=0}(p_j-1)^lS_{k-r-l}$.
 Thus
\beq
[\Db]=
(\prod p_i)\Bigg(\suml_{j=1}^r\frac{[X_j]}{p_j}\suml^{k-r}_{l=0}(p_j-1)^l\cS_{k-r-l}-c_1(\cN_{Z/M})\cS_{k-r}\Bigg),
\eeq
proving the first part of equation \eqref{Eq.Discriminant.Cohom.Class.Another.Form}.

\

The second part of equation \eqref{Eq.Discriminant.Cohom.Class.Another.Form} is now obtained immediately. For $\cS_{k-r}$ one uses
 property i. in \S\ref{Sec.Background.Combinatorics}, while  for  $\suml^{k-r}_{l=0}(p_j-1)^l\cS_{k-r-l}$ one uses
 property iv. in \S\ref{Sec.Background.Combinatorics}.
\epr

\subsection{An application: the multi-weighted-homogeneity of the classical discriminant of complete intersections}
\label{Sec.Proofs.Classical.Discriminant.Multi-Degree}
Consider an $r$-tuple of homogeneous forms in $k$ variables,
\beq
F_1(\ux)=\suml_{m_1+\cdots+m_k=p_1}a^{(1)}_{m_1\dots m_k}x^{m_1}_1\cdots x^{m_k}_k,
\quad \dots \quad
 F_r(\ux)=\suml_{m_1+\cdots+m_k=p_r}a^{(r)}_{m_1\dots m_k}x^{m_1}_1\cdots x^{m_k}_k.
\eeq
The classical discriminant of complete intersections is a polynomial in all these coefficients, $\frD\big(\{a^{(j)}_{m_1\dots m_k}\}_{j,\{m_I\}}\big)$.

The (non-)smoothness of the variety $\capl_j F^{-1}_j(0)\sset \P^{k-1}$ is obviously preserved under the scaling $\{F_j\to \la_j F_j\}$, for some
 non-zero constants $\{\la_j\in\k^*\}$. This induces the scaling of the coefficients,
  $\{a^{(j)}_{m_1\dots m_k}\}_{j,\{m_I\}}\to \{\la_j a^{(j)}_{m_1\dots m_k}\}_{j,\{m_I\}}$.

 Similarly, the (non-)smoothness is preserved under the scaling $x_i\to \la_i x_i$.
  It induces the scaling of the coefficients: $\{a^{(j)}_{m_1\dots m_k}\}_{j,\{m_I\}}\to \{\la^{m_1}_1\cdots\la^{m_k}_k\cdot a^{(j)}_{m_1\dots m_k}\}_{j,\{m_I\}}$.

Therefore the polynomial $\frD\big(\{a^{(j)}_{m_1\dots m_k}\}_{j,\{m_I\}}\big)$ is multi-weighted-homogeneous:
\beq\label{Eq.Discriminant.Classical.Weighted.Degrees.Benoist}\ber
\frD\big(\{\la_j a^{(j)}_{m_1\dots m_k}\}_{j,\{m_I\}}\big)=(\prodl^r_{j=1}\la_j^{deg_j})\cdot
\frD\big(\{a^{(j)}_{m_1\dots m_k}\}_{j,\{m_I\}}\big),
\\
\frD\big(\{\la^{m_1}_1\cdots\la^{m_k}_k\cdot a^{(j)}_{m_1\dots m_k}\}_{j,\{m_I\}}\big)=(\prodl^k_{i=1}\la_i)^{deg_{var}}
\frD\big(\{a^{(j)}_{m_1\dots m_k}\}_{j,\{m_I\}}\big).
\eer\eeq
(In the last formula one first gets the factor $\prod\la_i^{deg^{(i)}_{var}}$, by homogeneity. Then one gets $deg^{1}_{var}=\cdots=deg^{k}_{var}$
 as $\frD$ is invariant under the permutations of $x_1,\dots,x_k$ and the induced action on $\{a^{(j)}_{m_1\dots m_k}\}_{j,\{m_I\}}$.)

The degrees of this multi-weighted-homogeneity, $\{deg_j\}$, $deg_{var}$ are classically known in the hypersurface case ($r=1$).
 For complete intersections they have been recently computed, \cite[Theorem 1.3]{Benoist}.
 We recompute them as a simple application of the expression of the class $[\Db]$.

\bthe\label{Thm.Multidegrees.Classical.Discriminant}
1. $deg_j=\frac{(\prodl^r_{i=1}p_i)}{p_j}\suml^r_{i=1}\frac{(p_j-1)^k-(p_i-1)^k}{(p_j-p_i)}\frac{1}{\prodl_{l\neq i}(p_i-p_l)}$, for $j=1,\dots,r$.
\\2. $deg_{var}=\frac{1}{k}\suml^r_{j=1}p_jdeg_j=(\prodl^r_{j=1}p_j)\cdot\suml^r_{j=1}\frac{(p_j-1)^{k-1}}{\prodl_{l\neq j}(p_j-p_l)}$.
\ethe
Here the expression $\frac{(p_j-1)^k-(p_i-1)^k}{(p_j-p_i)}$ should be taken as $k(p_j-1)^{k-1}$ in the case $i=j$.
 The poles that seemingly occur when $p_i=p_j$ are canceled by grouping the terms
 $\frac{(p_j-1)^{k-1}}{\prodl_{l\neq j}(p_j-p_l)}+\frac{(p_i-1)^{k-1}}{\prodl_{l\neq i}(p_i-p_l)}$.
\bpr
{\bf 1.} For the sake of exposition we begin with the {\em hypersurface case}.

In this case the discriminant is the hypersurface in projective space, $\De=\{\frD(\{a_{m_1\dots m_k}\}_{\{m_I\}})=0\}\sset|\cO_{\P^{k-1}}(p)|$.
 The degree of $\De$ is the degree of the polynomial $\frD(\{a_{m_1\dots m_k}\}_{\{m_I\}})$, and it can be computed by intersecting $\De$ with
  the generic line $\P^1\sset |\cO_{\P^{k-1}}(p)|$. This line corresponds to a one-dimensional linear family of $deg=p$ hypersurfaces in $\P^{k-1}$.

We construct this one-dimensional family as follows.
 Let $Z=\P^1$ and $M=Z\times\A^{k}=\P^1\times\A^{k}$. Let $X\sset M$ be a hypersurface with $Sing(X)=Z\times\{0\}\sset M$.
 Assume $X$ is generically ordinary along $Z\times\{0\}$,    of generic multiplicity $p$.
  Assume  $\cO_M(X)\approx\cO_{\P^1}(1)\boxtimes\cO_{\A^k}$, i.e. the defining equation of $X\sset M$ is linear in homogeneous coordinates on $\P^1$.
  Then $\P\cN_{Z/X}=\tX\cap E\to Z$ is a linear family of projective hypersurfaces of degree $p$ in $\P^{k-1}$.
    The generic member of this family is smooth and the singular fibres, over $\Db\sset \P^1$, correspond to the intersection points of this line with
     the classical discriminant, $\De\sset |\cO_{\P^{k-1}}(p)|$.
 Thus the     discriminantal schemes are identified,
     $Z\supset\Db\approx\De|_{\P^1}$.

If we take $X$ generic for the given data, then $\tX\cap E$ is smooth. Thus $deg(\De)=deg(\Db)$ can be obtained from the
 equivalence class of theorem \ref{Thm.Discriminant.Cohomology.Class.via.Thom.Porteus}.
 In our case $\cO_Z(X)=\cO_M(X)|_Z\approx\cO_{\P^1}(1)$, thus $deg[X|_Z]=1$. Moreover, $c_1(\cN_{Z/M})=0$.
Therefore equation \eqref{Eq.Discrim.Class.for.Hypersurfaces} gives:  $deg(\De)=k(p-1)^{k-1}$.

\

{\em The case of complete intersections.}

Now the discriminant is a hypersurface in multi-projective space, $\De=\{\frD(\{a^{(j)}_{m_1\dots m_k}\}_{j,\{m_I\}})=0\}\sset\prod|\cO_{\P^{k-1}}(p_j)|$.
Then $\{deg_j\}$ is the multi-degree of $\De$. Equivalently, this gives the class $[\De]\in Pic\Big(\prod|\cO_{\P^{k-1}}(p_j)|\Big)\approx \Z^{\oplus r}$.

 Fix some $1\le j\le r$, we compute $deg_j$ as follows. Fix some points $\{q^{(i)}\in |\cO_{\P^{k-1}}(p_i)|\}_{i\neq j}$. This defines the embedding:
 \beq
|\cO_{\P^{k-1}}(p_j)|\stackrel{\nu_j}{\hookrightarrow}\prodl^r_{i=1}|\cO_{\P^{k-1}}(p_i)|,\quad
|\cO_{\P^{k-1}}(p_j)|\stackrel{\nu_j}{\hookrightarrow}\{q^{(1)}\}\times\cdots\times \{q^{(j-1)}\}\times |\cO_{\P^{k-1}}(p_j)|\times
\{q^{(j+1)}\}\times\cdots\times \{q^{(k)}\}.
\eeq
Then a line $\P^1\sset |\cO_{\P^{k-1}}(p_j)|$ is mapped to $\nu_j(\P^1)\sset \prodl^r_{i=1}|\cO_{\P^{k-1}}(p_i)|$.
Suppose both $\{q^{(i)}\in |\cO_{\P^{k-1}}(p_i)|\}_{i\neq j}$ and $\P^1\sset |\cO_{\P^{k-1}}(p_j)|$ are generic, then $\nu_j(\P^1)$ intersects $\De$
 transversally and $deg_j=deg(\nu_j(\P^1)\cap \De)$.

 This line corresponds to a one dimensional linear family of complete intersections,
\beq
\{f_1(\ux)=\cdots=f_{j-1}(\ux)=f_j(\ux,(t_0,t_1))=f_{j+1}(\ux)=\cdots=f_r(\ux)=0\}\sset \P^{k-1}\times\P^1.
\eeq
(Here $f_j(\ux,(t_0,t_1))$ is linear in $(t_0,t_1)$.)

As in the hypersurface case, we realize this line as the singular locus, as follows.
Let $Z=\P^1$ and $M=Z\times\A^k=\P^1\times\A^k$. Fix an $r$-tuple of hypersurfaces, $\{X_i\sset M\}_{i=1\dots r}$ satisfying:
\beq
\cO_M(X_i)=\Big\{\ber \cO_{\P^1}\boxtimes\cO_{\A^k},\ i\neq j\\\cO_{\P^1}(1)\boxtimes\cO_{\A^k},\ i=j\eer,\quad Sing(\capl^r_{i=1} X_i)=Z\times\{\o\}.
\eeq
We assume that the multiplicity of $X_j$ along $Z\approx Sing(\cap X_i)$ is constant.
(For the other $X_i$ this holds, as their defining equations do not depend on the coordinates of $Z$.)
 Moreover, we assume that $\capl^r_{i=1} X_i$ is s.c.i. at each point of $Z$ and is generically ordinary along $Z$.
Therefore the line $\nu_j(Z)\sset \prodl^r_{i=1}|\cO_{\P^{k-1}}(p_i)|$ does not lie inside $\De$, and $\Db\sset Z$ coincides with $\De|_{\nu_j(Z)}$.

 Finally, $deg(\Db)$ is computed by theorem \ref{Thm.Discriminant.Cohomology.Class.via.Thom.Porteus}. By our construction:
 \beq
 \cO_Z(X_i)=\Big\{\ber \cO_{\P^1},\ i\neq j\\ \cO_{\P^1}(1),\ i=j\eer,\quad c_1(\cN_{Z/M})=0.
\eeq
Therefore equation \eqref{Eq.Discriminant.Cohom.Class.Another.Form} gives:
\beq
deg_j=[\Db]=\frac{(\prod p_i)}{p_j}\suml^r_{i=1}\frac{(p_j-1)^k-(p_i-1)^k}{(p_j-p_i)\prodl_{l\neq i}(p_i-p_l)}.
\eeq

The reproves equation (5) of Theorem 1.3 in \cite{Benoist}.

{\bf 2.} We study a particular class of monomials in the polynomial $\frD$. Consider the intersection of hypersurfaces
\beq
\Big(\suml^k_{i=1}a_i^{(p_1)}x_i^{p_1}=0\Big)\cap
\Big(\suml^k_{i=1}a_i^{(p_2)}x_i^{p_2}=0\Big)\cap\dots\cap
\Big(\suml^k_{i=1}a_i^{(p_r)}x_i^{p_r}=0\Big)\sset \P^{k-1}.
\eeq
Suppose none of the coefficients $\{a_i^{(p_j)}\}$  vanishes, then this is the intersection of smooth hypersurfaces. Moreover, if the coefficients $\{a_i^{(p_j)}\}$
 are mutually generic, then this variety is a smooth projective complete intersection. Therefore, despite the fact that we have put ``most" of
  the coefficients  $\{a_{m_1..m_k}^{(p_j)}\}$ to zero, the restricted discriminantal polynomial $\frD(\{a_i^{(p_j)}\}_{i,j})$ does not vanish identically.

For every $j$ this polynomial is of degree $deg_j$ in the set $\{a^{(p_j)}_1,\dots,a^{(p_j)}_k\}$ of coefficients. Suppose a monomial
\beq
\Big(\prodl_i (a^{(p_1)}_i)^{n^{(1)}_i}\Big)\cdot \Big(\prodl_i (a^{(p_2)}_i)^{n^{(2)}_i}\Big)\cdots
\Big(\prodl_i (a^{(p_r)}_i)^{n^{(r)}_i}\Big)
\eeq
participates in this polynomial. Then, for any $j$ we have: $\suml^k_{i=1}n^{(j)}_i=deg_j$. Under the scaling $a^{(p_j)}_i\to \la^{p_j}_ia^{(p_j)}_i$
 this monomial is multiplied by $\prodl^r_{j=1}\prodl^k_{i=1}\la^{p_j n^{(j)}_i}_i$. And by the weighted homogeneity,
  equation \eqref{Eq.Discriminant.Classical.Weighted.Degrees.Benoist},
 this factor equals $(\prodl_i\la_i)^{deg_{var}}$. Checking the exponent of each $\la_i$ we get:
 $\suml^r_{j=1}p_j n^{(j)}_i=deg_{var}$, for each $i$. Finally, we sum over all $i$ to get:
\begin{multline}
k\cdot deg_{var}=\suml^k_{i=1}\suml^r_{j=1}p_j n^{(j)}_i=\suml^r_{j=1}(\suml^k_{i=1}p_j n^{(j)}_i)=\suml^r_{j=1}p_jdeg_j
\stackrel{part \ 1}{=\!=}\\
=
\suml^r_{j=1}\Big(\prodl^r_{i=1}p_i\Big)\suml^r_{i=1}\frac{(p_j-1)^k-(p_i-1)^k}{(p_j-p_i)}\frac{1}{\prodl_{l\neq i}(p_i-p_l)}
\stackrel{\substack{Part\ \rm{vi.}\ in \ \S\ref{Sec.Background.Combinatorics}}}{=\!=\!=\!=}
(\prodl^r_{j=1}p_j)\cdot\suml^r_{j=1}\frac{k(p_j-1)^{k-1}}{\prodl_{l\neq j}(p_j-p_l)}.\quad
\end{multline}
This reproves equation (6) of Theorem 1.3 in \cite{Benoist}.\epr

\subsubsection{An alternative derivation in the hypersurface case} Consider the hypersurface case, i.e. the complete linear system
 $|\cO_{\P^k}(p)|$. Of the degrees $\{deg_j\}$ we have only one, it is the degree of the discriminant, known by the classical
 formula of Boole \cite[pg. 38]{Gelfand-Kapranov-Zelevinsky}: $deg=k(p-1)^{k-1}$. The degree $deg_{var}$ can be also computed in various ``pedestrian" ways. For completeness we give one of them.

The number $deg_{var}$ is determined if we know at least one monomial of the polynomial $\frD(\{a_{m_1,..,m_k}\})$.
 In particular, we would like to show the presence of a monomial of the form: $\Big(a_{p,0,..,0}a_{0,p,0,..,0}\cdots a_{0,..,0,p}\Big)^N$, for some $N$.
  The degree  is fixed by the total degree of the discriminant, therefore $N=(p-1)^{k-1}$.

\

The monomial $\Big(a_{p,0,..,0}a_{0,p,0,..,0}\cdots a_{0,..,0,p}\Big)^{(p-1)^{k-1}}$
always participates in the discriminant. This follows immediately from Theorem 3.2 on pg. 362 of
\cite{Gelfand-Kapranov-Zelevinsky},
this monomial corresponds to a trivial triangulation of the polytope $Conv(x^p_1,\dots,x^p_k)$,
consisting just of the polytope itself.

Alternatively, the presence of this monomial can be shown directly as follows. Note that the hypersurface $\{\sum^k_{i=1} x^p_i=0\}\sset\P^{k-1}$
 is
smooth and for it the only non-vanishing monomials are of the type $a_{0,\dots,0,p,0,\dots,0}$.
Therefore in the polynomial that defines the discriminant there must be a monomial that contains only
 the coefficients $\{a_{p,0,\dots,0,},a_{0,\dots,0,p,0,\dots,0},a_{0,\dots,0,p}\}$ and no other $a_{*}$'s. Further, as the
 hypersurface  $\{\sum^k_{i=2} x^p_i=0\}$ is singular, the monomial above must contain {\em all} of
 $\{a_{0,\dots,0,p,0,\dots,0}\}$. The permutation group $\Si_k$
acts on the coordinates, therefore this monomial must be of the form:
$\Big(a_{p,0,..,0}a_{0,p,0,..,0}\cdots a_{0,..,0,p}\Big)^{N}$. Finally, we get, as before
$N=(p-1)^{k-1}$.

Therefore we get: $deg_{var}=p(p-1)^{k-1}$.

\subsection{Proof of Theorem \ref{Thm.Discriminant.Cohomology.Class}  (the general case)}
\label{Sec.Proofs.Cohomology.Class.Discriminant.Classical.Style}
Now we compute the class $[\Db]\in\Pic(Z)$ via the weighted homogeneity of the classical discriminant.

\bpr
{\bf Step 1.}
By the assumption $X$ is a globally complete intersection, $X=\capl_j X_j\sset M$. Take a global section  $f_j\in H^0(\cO_M(X_j))$.

Take a point $\o\in Z$, take the locally complete intersection presentation, $(Z,\o)=\capl_i(Z_i,\o)\sset(M,\o)$.
  Let $g_i\in H^0(\cO_{(M,\o)}(Z_i,\o)$ be a section that locally defines $Z_i$.

By the assumption, the generic multiplicity of $X_j$ along $Z$ is $p_j$, thus $f_j\in I^{(p_j)}_{(Z,\o)}$, the symbolic power of ideal,
 \cite[\S2.2]{Kerner.Nemethi.Discriminant.Local}. As $(Z,\o)$ is a complete intersection, the symbolic and ordinary powers coincide,
 $I^{(p_j)}_{(Z,\o)}=I^{p_j}_{(Z,\o)}$, see e.g. lemma 2.3 of \cite{Kerner.Nemethi.Discriminant.Local}.
Thus  we have the expansion
\beq
f_j=\suml_{m_1+\cdots+m_k=p_j} g^{m_1}_1\cdots g^{m_k}_k a^{(j)}_{m_1,..,m_k},\quad
  \{a^{(j)}_{m_1..m_k}\in \cO_{(M,\o)}\}_{j,I}.
\eeq
The coefficients $\{a^{(j)}_{m_1,\dots,m_k}\}$ are not defined uniquely, because of  Koszul relations among $\{g_i\}$. But
 their restrictions, $\Big\{a^{(j)}_{m_1,\dots,m_k}|_{(Z,\o)}\in\cO_{(Z,\o)}\Big\}$, are uniquely defined, see \S\ref{Sec.Background.Normal.Cone}.

\

By the choice of $\{f_j\}$, $\{g_i\}$ the coefficient $a^{(j)}_{m_1,..,m_k}|_{(Z,\o)}$ is a local section of the line
 bundle $\cO_{(Z,\o)}(X_j-\sum m_i Z_i)$.

The full collection $\{a^{(j)}_{m_I}|_{(Z,\o)}\}_{j,I}$ then glues to a (local) section of the bundle
\begin{multline}
\oplusl_{j,\{m_I\}}\cO_{(Z,\o)}(X_j-\sum m_i Z_i)=Hom\Big(Sym^p\big(\oplus \cO_{(Z,\o)}( Z_i)\big),\oplusl_j\cO_{(Z,\o)}(X_j)\Big)=\\=
Hom\Big(Sym^p(\cN_{Z/M})|_{(Z,\o)},\cN_{X/M}|_{(Z,\o)}\Big).
\end{multline}

We claim that as the point $\o$ varies over $Z$ these local sections glue to a global section of
\beq
Hom\Big(Sym^p(\cN_{Z/M}),\cN_{X/M}|_Z\Big).
\eeq
This happens because all the constructions are functorial.
 Indeed, fix two points $\o,\o'\in Z$, we should check the transition maps from $\{a^{(j)}_{m_I}|_{(Z,\o)}\}_{j,I}$ to
  $\{a^{(j')}_{m_{I'}}|_{(Z,\o')}\}_{j',I'}$. Fix some common neighborhood $\o,\o'\in\cU\sset Z$, such that both $\{g_i\}$ and $\{g'_{i'}\}$ are defined on $\cU$.
   Then the transition from $\{g_i\}$ to $\{g'_{i'}\}$ is done by an element $\phi_{\o,\o'}\in GL_k(\cU)$.
    The elements $\{\phi_{\o,\o'}\}_{\o,\o'\in Z}$ are the transition maps for $\cN_{Z/M}$.
    Similarly one trivializes $\cN_{X/M}|_Z$ to get the transition maps $\{\psi_{\o,\o'}\}_{\o,\o'\in Z}$.
 Accordingly, the transition from  $\{a^{(j')}_{m'_1,..,m'_k}\}$ to $\{a^{(j)}_{m_1,..,m_k}\}$ is done by $(\phi^{\otimes p}_{\o,\o'})^\vee\otimes \psi_{\o,\o'}$.
 But these are precisely the transition maps for the bundle $Hom\Big(Sym^p(\cN_{Z/M}),\cN_{X/M}|_Z\Big)$.

\

Finally, we compose this (global) section with the discriminantal map of \S\ref{Sec.Proofs.Classical.Discriminant.Multi-Degree}:
\beq
\prodl^r_{j=1} H^0(\cO_{\P^k}(p_j))\stackrel{\frD}{\to}\k,\quad\quad
\{a^{(j)}_{m_1\dots m_k}|_Z\}_{j,\{m_I\}}\to \frD\Big(\{a^{(j)}_{m_1\dots m_k}|_Z\}_{j,\{m_I\}}\Big).
\eeq
It remains to verify that $\frD(\{a^{(j)}_{\{m_I\}}|_Z\}_{j,\{m_I\}})$ is
 a section of some well defined line bundle on $Z$, and to identify this bundle. For this we should check that all
  the monomials of this polynomial are sections of the same line bundle. More precisely, we should check that every monomial
   of   $\frD(\{a^{(j)}_{\{m_i\}}|_Z\}_{j,\{m_i\}})$ is a section of some $\cO_Z(\sum a_j X_j-\sum b_i Z_i)$, where the
    constants $\{a_j\}$, $\{b_i\}$ are the same for all the monomials, and to find these constants.

{\bf Step 2.}
 As is explained
 in \S\ref{Sec.Proofs.Classical.Discriminant.Multi-Degree}, $\frD$ is multi-weighted-homogeneous,  see
 equation \eqref{Eq.Discriminant.Classical.Weighted.Degrees.Benoist}.
 Thus we get immediately: $a_j$, $b_i$ do not depend on a particular monomial and
\beq
a_j=deg_j,\quad b_i=deg_{var}\ \text{(independently of $i$)}.
\eeq
The degrees $deg_j$, $deg_{var}$ are computed in theorem \ref{Thm.Multidegrees.Classical.Discriminant}. Thus we get:
  if locally $(Z,\o)=\cap(Z_i,\o)$ then
$\frD\Big(\{a_{m_1,..,m_k}|_{(Z,\o)}\}_{\substack{\{m_j\}}}\Big)$ is the section of the bundle
\beq
\cO_{(Z,\o)}\Bigg((\prodl^r_{i=1}p_i)\suml^r_{j=1}\Big(\frac{X_j}{p_j}\suml^r_{i=1}\frac{(p_j-1)^k-(p_i-1)^k}{(p_j-p_i)\prodl_{l\neq i}(p_i-p_l)}
-\frac{(p_j-1)^{k-1}}{\prodl_{i\neq j}(p_j-p_i)}(\suml^k_{l=1}Z_l)\Big)\Bigg).
\eeq
Globally these sections glue to a section of
\beq
\cO_Z\Bigg((\prodl^r_{i=1}p_i)\suml^r_{j=1}\Big(\frac{X_j}{p_j}\suml^r_{i=1}\frac{(p_j-1)^k-(p_i-1)^k}{(p_j-p_i)\prodl_{l\neq i}(p_i-p_l)}
-\frac{(p_j-1)^{k-1}}{\prodl_{i\neq j}(p_j-p_i)}det(\cN_{Z/M})\Big)\Bigg).
\eeq
This proves the statement.
\epr

\subsection{The case: $X\sset M$ is a locally complete intersection}\label{Sec.X.locally.Complete.Intersection}
Here we provide the details to the remark \ref{Remarks}  and extend the formula of theorem \ref{Thm.Discriminant.Cohomology.Class}
to this case, i.e. prove the formula of equation \eqref{Eq.Discriminant.Cohom.Class.loc.complete.inters}.

Recall (\S\ref{Sec.Background.Discriminant}) that given a subscheme $Z\sset M$ one passes from the sheaf $\cO_M$ to the associated graded sheaf,
\[
gr_Z\cO_M=\quotients{\cO_M}{I_Z}\oplus\quotients{I_Z}{I^2_Z}\oplus\cdots.
\]
As $Z$ is a locally complete intersection, the stalks of this sheaf are graded algebras, $(gr_Z\cO_M)|_{(Z,\o)}\approx\cO_{(Z,\o)}[\uy]$.
 In particular, the restriction/transition
 morphisms are graded. For a subsheaf $I_X\sset\cO_M$ we get the natural graded subsheaf $gr_Z(I_X)\sseteq gr_Z(\cO_M)$.

\

 Let $Z\sseteq Sing(X)$ be a connected component and suppose that for each point $\o\in Z$ the germ $(X,\o)$ is s.c.i. over $(Z,\o)$, with the
 multiplicity sequence $p_1\le\dots\le p_r$. By the connectedness and the s.c.i. condition this multiplicity sequence  is the
  same on all the components of  $Z$, see Example \ref{Ex.SCI.hypersurface}.
   We construct the filtration of $gr_Z(I_X)$ by
   graded subsheaves, $J_{p_1}\sseteq J_{p_2}\sseteq\cdots\sseteq J_{p_r}=gr_Z(I_X)$, as follows.

For any $o\in Z$ fix a good basis of $I_{(X,\o)}$, as in \S\ref{Sec.Background.s.c.i}, and take its
 graded image, i.e. the leading terms of the generators, $gr_{(Z,\o)}(I_{(X,\o)})=(\tf_{p_1},\dots,\tf_{p_r})\sset \cO_{(Z,\o)}[\uy]$.
Suppose $p_1=\cdots=p_{k_1}<p_{k_1+1}$ and define
\beq
J^{(Z,\o)}_{p_1}:=J^{(Z,\o)}_{p_2}:=\cdots=J^{(Z,\o)}_{p_{k_1}}:=(\tf_{p_1},\dots,\tf_{p_{k_1}})\sset\cO_{(Z,\o)}[\uy].
\eeq
 This is a graded submodule. It does not depend on the choice of generators as $\tf_{p_1},\dots,\tf_{p_{k_1}}$ cannot be mixed with $\tf_{p_j}$ for $j>k_1$.

 We claim that the modules $\{J^{(Z,\o)}_{p_1}\}_{\o\in Z}$ glue to a graded subsheaf of modules, $J_{p_1}\sseteq gr_Z I_X$.
It is enough to check the transitions/restrictions of the generators of $J^{(Z,\o)}_{p_1}$. For any $\o,\o'\in Z$ the embedding
 $(Z,\o)\cap(Z,\o')\into (Z,\o)$ induces the pullback $gr_Z I_X|_{(Z,\o)}\stackrel{i^*}{\to} gr_Z I_X|_{(Z,\o)\cap(Z,\o')}$.
 We should verify that $i^*$ sends  $J^{(Z,\o)}_{p_1}$ to  $J^{(Z,\o)\cap(Z,\o')}_{p_1}$, without mixing it with the other elements of $gr_Z(I_X)$.
 Indeed, the sheaf $gr_Z I_X$ is graded, thus
 $i^*(\tf_{p_1}),\dots,i^*(\tf_{p_{k_1}})$ are still homogeneous, of degree $p_1$. Therefore,
  being of the lowest degree, these must lie in $J^{(Z,\o)\cap(Z,\o')}_{p_1}$.

\

In the same way, if $p_{k_1+1}=\cdots=p_{k_2}<p_{k_2+1}$, one defines the subsheaf $J_{p_{k_1+1}}=\cdots=J_{p_{k_2}}=(\tf_{p_1},\dots,\tf_{p_{k_2}})\sset gr_Z(I_X)$, and so on.
 As before, each of these subsheaves does not depend on the choice of generators.
This defines the global filtration of $gr_Z I_X$ by sheaves of ideals,
\beq
0= J_{p_0}\subsetneq J_{p_1}\sseteq J_{p_2}\sseteq\cdots\sseteq J_{p_r}=gr_Z(I_X).
\eeq
Thus we have a filtration of $\quotients{gr_Z(I_X)}{gr_Z(I_X)^2}$ by {\em locally free} sheaves:
\beq
0\sseteq\quotient{J_{p_1}+gr_Z(I_X)^2}{gr_Z(I_X)^2}\sseteq\cdots \quotient{J_{p_i}+gr_Z(I_X)^2}{gr_Z(I_X)^2}\sseteq\cdots\sseteq \quotient{gr_Z(I_X)}{gr_Z(I_X)^2}.
\eeq

To define the  filtration of the restricted normal bundle, $\cN_{X/M}|_Z$, we use the identification
 $\cN_{X/M}|_Z\approx (gr_Z \quotients{I_X}{I^2_X})^\vee=Hom(gr_Z\quotients{I_X}{I^2_X},gr_Z\cO_X)$.
  Thus we get:
\beq
\cN_{X/M}|_Z=\cN_{p_0}\supsetneq\cN_{p_1}\supseteq\cN_{p_2}\supseteq\cdots\supseteq\cN_{p_r}=0,\quad
 \text{  where }
  \cN_{p_j}=\{s\in Hom(gr_Z\quotients{I_X}{I^2_X},gr_Z\cO_X)|\ s(J_{p_j})=0\}.
\eeq
 Accordingly,  $det(\cN_{X/M}|_Z)=\otimesl^r_{i=1} det\quotients{\cN_{p_{i-1}}}{\cN_{p_i}}$.

In the case of a globally complete intersection,
 $(X,Z)=\capl(X_j,Z)$, the normal bundle splits:  $\cN_{X/M}|_Z\approx\oplus \cO_Z(X_j)$.
 Here the filtration is: $\cN_{p_i}=\oplusl_{p_j> p_i}\cO_Z(X_j)$.

 Thus, by the splitting principle, it is enough  to replace in the formula of theorem \ref{Thm.Discriminant.Cohomology.Class}
  the part $\suml^r_{j=1}\frac{[X_j]}{p_j(1-p_jt)}$
  by $\suml^r_{j=1}\frac{[det(\quotients{\cN_{p_{j-1}}}{\cN_{p_{j}}})]}{p_j(1-p_jt)}$.

\section{Examples, classes of some strata of $\Db$ and applications}\label{Sec.Examples.Applications}

\subsection{Examples}
In the hypersurface case, $r=1$, we get: $\cO_Z(\Db)=\cO_Z\Big((p-1)^{k-1}\big(kX-p\cdot det(\cN_{Z/M})\big)\Big)$.
 We consider the following particular cases:
\bee[i.]
\item Let $Z\sset\P^{k+1}$ be a projective curve, a
reduced complete intersection of multidegree $(d_1,\ldots,d_k)$. Suppose
 $Z=Sing(X)$, where
$X\sset\P^{k+1}$ is a hypersurface of degree $d$. Then
\[
{\rm deg}[\Db]=(\prod d_i)(p-1)^{k-1}(kd-p\sum d_i).
\]

\item Let $M$ be the blowup of $\P^{k+1}$ at a point $p$, with the exceptional divisor $E\approx
\P^k$. Let $Z\sset E$ be a smooth  complete
intersection of dimension one  of multidegree $(d_1,\ldots
,d_{k-1})$. Suppose $X\sset M$ is an irreducible hypersurface, $Z=Sing(X)$  and
the intersection $X\cap E$ is of degree $q$ in $E$. Then
\[
\deg [\Db]=(\prod d_i)(p-1)^{k-1}(kq+p-p\sum d_i).
\]
\item  Let $X=\cupl^p_{i=1}X_i\sset M$ be the union of smooth
hypersurfaces intersecting pairwise-transversally. Suppose that
for any pair $i\neq j$ the intersection  $X_i\cap X_j$ is the same
subset of $M$, denote it by $Z$. Hence $Z$ is smooth and the transversal
type of $X$ along $Z$, at any point, is an \omp\
of multiplicity $p$. Thus we expect: $\Db=\empty$.
Compare this expectation with the class $[\Db]$ obtained from
Theorem \ref{Thm.Discriminant.Cohomology.Class}. We have ${\rm codim}_M Z=2$,
 $\det(\cN_{Z/M})=[X_1]+[X_2]$, thus
\[[\Db]=(p-1)\Big(2\sum^p_{i=1}[X_i]-p([X_1]+[X_2])\Big)=0.
\]
 (For $p>2$ we use the isomorphism $\cO_Z(X_i)\approx\cO_Z(X_j)$, which is proved as follows.
As all the intersections $X_i\cap X_j$ are transversal, we have:
$\cO_Z(X_i)\oplus\cO_Z(X_k)\approx\cN_{Z/M}\approx \cO_Z(X_j)\oplus\cO_Z(X_k)$, for any $j\neq  k\neq i$. This is an isomorphism of bundles of rank two.
 Then the second exterior powers are isomorphic,
 $\cO_Z(X_k+X_i)\approx \cO_Z(X_k+X_j)$. Tensor by $\cO_Z(-X_k)$ to get $\cO_Z(X_i)\approx \cO_Z(X_j)$.)

 As a non-empty discriminant must be an effective divisor (and then $[\Db]\neq0$), we get again: $\Db=\empty$.
 \eee

\subsection{Classes of some further strata}\label{Sec.Cohomol.Classes.of.Further.Strata}
Let $X\sset M$ be a hypersurface and $Z=Sing(X)$. Suppose $dim(Z)>1$, so the discriminant $\Db$ is stratified, \S\ref{Sec.Introduction.Equivalence.Classes.Strata}.
In addition to the assumptions of \S\ref{Sec.Introduction.Assumptions},
 suppose the  strata $\Si_{A_2}$, $\Si_{A_1A_1}$ are non-empty and of codimension two in $Z$. We compute their classes in the Chow group $A^2(Z)$.
\bthe
Set $c_i=c_i(\cN_{Z/M})$. Then the classes of (the closures of) these strata, in $A^2(Z)$,  are given by:
\beq
[\lSi_{A_2}]=\frac{(k-1)(p-2)(p-1)^{k-2}}{2}\Big([c_1]^2p(2p-1)+[c_2]2p(\frac{2k}{k-1}-p)-2[c_1][X]p(k+1)+[X]^2k(k+1)\Big)
\eeq
\begin{multline}
[\lSi_{A_1,A_1}]=\frac{(p-1)^{2k-2}}{2}[Xk-c_1p]^2+
\\
+\frac{(p-1)^{k-2}}{4}\Bigg(\hspace{-0.2cm}
\ber [X]^2k\big(6k^2-4-p(k+1)(3k-2)\big)+2p[X][c_1]\big(p(k+1)(3k-2)-2(3k^2-2)\big)+
\\+2p[c_2]\big(p^2(3k-2)-4p(3k-1)+12k\big)+2p[c_1]^2\big(p(9k-8)-6(k-1)-p^2(3k-2)\big)\eer
\hspace{-0.3cm}\Bigg).
\end{multline}
\ethe
\bpr
As in the proof of theorem \ref{Thm.Discriminant.Cohomology.Class} we use the diagram
\beq\label{Eq.Diagram.for.Critical.Locus.Extended}
\bM Crit(\pi)\sset \tX\cap E\stackrel{\nu}{\hookrightarrow}E&\approx \P\cN_{Z/M}&\sset Bl_Z(M)\\
\hspace{2.2cm}\pi\searrow\hspace{0.3cm}
\downarrow&&\downarrow\\\hspace{3cm}Z&\sset& M
\eM
\eeq
 By the general theory,  the classes of the
 strata are obtained by specializing the relevant Thom polynomials (e.g. see Table 1 of \cite{Kazarian2000}
 and  Table 2 of \cite{Kazarian2006}). In particular, in the low co-dimension cases we have the following
 classes in $A^*(Crit(\pi))$:
\beq\ber
[\lSi_{A_1}]=1,\quad [\lSi_{A_2}]=a_1,\quad [\lSi_{A_1,A_1}]=\frac{c_1(\cN_{\De^\bot/Z})-u-3a_1}{2},\quad
  [\lSi_{A_3}]=ua_1+3a_2,\quad [\lSi_{D_4}]=a_1a_2-ua_2-2a_3.
\eer\eeq
Here $u=[\pi^*X]-p[E]$, while the classes $\{a_i\}$ are defined by the formula:
\beq
\suml_i a_i=c(T_{E/Z}^*\otimes\cO_E(\tX)-T_{E/Z})=
\frac{(1+u+E)^k-(1+u+E)^{k-1}c_1+\cdots\pm c_k}{\Big((1-E)^k+(1-E)^{k-1}c_1+\cdots+c_k\Big)}\quad\text{ for }
c_i=c_i(\cN_{Z/M}).
\eeq
As in \S\ref{Sec.Proofs.Cohomology.Class.Discriminant.Thom.Porteus}  we pushforward these classes to $A^*(E)$, i.e.
multiply these expressions by $[Crit(\pi)]\in A^1(\tX\cap E)$ and by $[\tX\cap E]=[\pi^*X]-p[E]$.
Note that the class $[Crit(\pi)]$ in the hypersurface case is
\beq
c_k\Big(\pi^*\cN^*_{Z/M}\otimes\cO_E(\tX)(-1)\Big)=(u+E)^k-(u+E)^{k-1}c_1(\cN_{Z/M})+\cdots\pm c_k(\cN_{Z/M}).
\eeq
Then one applies the Gysin homomorphism to get the classes in $A^*(Z)$.
\epr

\bex
Similarly to \S\ref{Sec.Proofs.Classical.Discriminant.Multi-Degree}, in the simplest case we recover the degrees of the strata
  of classical discriminant. Let $Z=\P^2$ and $M=Z\times\A^k$. Suppose $\cO_M(X)=\cO_{\P^2}(1)\boxtimes\cO_{\A^k}$ and $Sing(X)=Z\times\{\o\}\sset Z\times\A^k$.
   Then $deg[X]=1$ and $det(\cN_{Z/M})=0$ and we have:
\bee[i.]
\item
$deg(\lSi_{A_2})=\frac{k(k+1)(k-1)(p-2)(p-1)^{k-2}}{2}$, which recovers the formula (``$deg(C)=\dots$") on page 10 of \cite{Aluffi-1998},
 see also equation (60) of \cite{Kerner.2008}.
\item
$deg(\lSi_{A_1,A_1})=\frac{k^2(p-1)^{2k-2}}{2}+\frac{k(p-1)^{k-2}}{4}\Bigg( 6k^2-4-p(k+1)(3k-2)\Bigg)$, which recovers
 the formula (``$deg(G)=\dots$") on page 10 of \cite{Aluffi-1998}.
\eee
\eex

\subsection{A bound on jumps of multiplicity}\label{Sec.Application.Bound.on.jump.of.multiplicity}
Let $X\sset M$ be a hypersurface of dimension $n$.
\marginpar{{\bf WHY not k instead of k???}}
Suppose $Z=Sing(X)$  is a smooth projective irreducible curve, the generic
multiplicity of $X$ along $Z$ is $p$ and the generic transversal type is ordinary.
\bprop
Suppose that at some points of $Z$ the multiplicity jumps, $\{mult(X,\p_i)=p+p_i\}_i$. Then
\[\suml_i p_i\le deg([X]|_Z)- \frac{p}{n}deg(\det(\cN_{Z/M})).
\]
\eprop
\bpr
The discriminant $\De^\bot\sset Z$ is zero dimensional, thus it is enough to compute the local degrees of $\De^\bot$ at
the points $\{\p_i\}_i$.

Choose the local coordinates  such that
$(Z,\p_i)=\{x_1=\cdots=x_n=0\}\sset(M,\p_i)$ and $x_{n+1}$ is the local coordinate along $(Z,\p_i)$.
(If needed, we can pass to completion,  $\widehat{(Z,\p_i)}$.)
As $mult(X,\p_i)=p+p_i$, this hypersurface singularity can be presented in the form
\beq\label{Eq.Multiplcity.Jumps.Bound}
\{x^{p_i}_{n+1}f_p(x_1,\dots,x_n)+g(x_1,\dots,x_n,x_{n+1})=0\}\sset(M,\p_i),
\eeq
where $f_p$ is a homogeneous form of degree $p$, while
\begin{multline}
g(x_1,\dots,x_n,x_{n+1})\in (x_1\dots x_n)^{p+1}\cap\cm^{p+p_i}+(x_1\dots x_n)^{p}\cap(x_{n+1})^{p_i+1}=\\=
(x_1\dots x_n)^{p+1}\cap\cm^{p+p_i}+(x_1\dots x_n)^{p}\cap\cm^{p+p_i+1}.
\end{multline}
Deform the function in this equation (locally), to assume that $f_p$ is non-degenerate. This induces a flat deformation of $\Db$.
 So the local degree of $\De^\bot(X,\p_i)$
 is at least the local degree of the discriminant for this particular deformed singularity. We should compute the degree of the discriminant for the projection
\beq
(Z,\o)\times \P^{n-1}\supset\tX\cap E=\{x^{p_i}_{n+1}f_p(x_1,\dots,x_n)=0\}\to (Z,\o)=Spec(\k[[x_{n+1}]]).
\eeq
This computation can be done in either of the following ways.
\bee[i.]
\item
The needed degree is $p_i$ times the degree for the projection $\tX\cap E=\{x_{n+1}f_p(x_1,\dots,x_n)=0\}\to (Z,\o)$.
(This can be seen e.g. by deformation $\{(x^{p_i}_{n+1}-\ep)f_p(x_1,\dots,x_n)=0\}$.) For the later (reduced) case we notice that the only
singularity occurs over $\o\in Z$. Thus the curve $Z$ is not important (neither locally nor globally), so we replace $Z$ by $\P^1$.
 Now we have a linear family of projective hypersurfaces in $\P^{n-1}$, of degree $p$. For $x_{n+1}\neq 0$ all the hypersurfaces
 are smooth (and of the right dimension). This family corresponds to a line in the parameter space of the hypersurfaces.
  And the local degree of $\Db$, at $x_{n+1}=0$, equals the total intersection number of this line with the
   classical discriminant. But the later number is just $n(p-1)^{n-1}$. Thus the local degree of $\Db$ at $\o\in Z$ equals $p_in(p-1)^{n-1}$.
\item
As the form $f_p(x_1,\dots,x_n)$ is non-degenerate, the hypersurface singularity $\{f_p(x_1,\dots,x_n)=0\}\sset(\k^n,\o)$ is $\mu=const$ equisingular to
the Brieskorn-type singularity, $\{\suml^n_{j=1}a_j x^p_j=0\}$, with $\{a_j\neq0\}_j$. Thus we can assume $f_p$ in this form, and
 equation \eqref{Eq.Multiplcity.Jumps.Bound} becomes:
\beq
\suml^n_{j=1}x^{p_i}_{n+1}a_j x^p_j+g(x_1,\dots,x_n,x_{n+1})=0.
\eeq
Now deform this equation to: $\suml^n_{j=1}(x^{p_i}_{n+1}-\ep_j)a_j x^p_j+g(x_1,\dots,x_n,x_{n+1})=0$, where $\{\ep_j\}$ are all distinct.
This induces a flat deformation of the discriminant, splitting the initial
point into $n\cdot p_i$ discriminantal points.
 For each such point
the multiplicity of $X$ is locally constant along $Z$, thus
we can compute the local degree of the discriminant by using theorem \S 5.5 of \cite{Kerner.Nemethi.Discriminant.Local}.
Namely, it is the degree of the scheme
\beq
\{\di_1 f_p=\cdots=\di_n f_p=0\}=\{x^{p-1}_1=\cdots=x^{p-1}_{n-1}=tx^{p-1}_n=0\}\sset \P^{n-1}\times(Z,\o).
\eeq
This scheme is supported at one point, $(0,\dots,0,1)$, its degree is $(p-1)^{n-1}$. Thus the local degree of $\Db$ at $\o\in Z$ equals $p_in(p-1)^{n-1}$.
\eee

By going over all the discriminantal points on $Z$, we get:
$deg(\De^\bot)\ge  n(p-1)^{n-1}\suml_i p_i$.

Note that these local contributions do not depend on the global classes $\cN_{Z/M}$, $\cN_{X/M}$.

Now, comparing with equation \eqref{Eq.Discrim.Class.for.Hypersurfaces}, with $k=n$, we get:
\beq
\suml_i p_i\le deg([X]|_Z)- \frac{p}{n}deg(\det(\cN_{Z/M})).\quad \text{\epr}
\eeq

\end{document}